\documentclass[12pt, oneside]{article}
\usepackage{amsmath,amsfonts,amssymb,amscd,theorem}
\date{}

\newtheorem{thm}{Theorem}
\newtheorem{lem}{Lemma}
\newtheorem{defi}{Definition}
\newtheorem{cor}{Corollary}
\newtheorem{conj}{Conjecture}
\newtheorem{prop}{Proposition}

\begin{document}

\title{Cohomologies of unipotent harmonic bundles over quasi-projective
varieties I: \\The case of noncompact curves}
\author{J\"urgen Jost, Yi-Hu Yang\thanks{Supported partially
by NSF of China (No. 10471105) and "Shuguang Project" of Committee
of Education of Shanghai (04SG21)}, and Kang Zuo} \maketitle

\section{Introduction}
Let $\overline S$ be a compact Riemann surface (holomorphic curve) of genus $g$. Let $p_1,
p_2, \cdots, p_s$ be $s>0$ points on it; these points define a divisor, and we denote the open Riemann surface
${\overline S}\setminus\{p_1, \dots, p_s\}$ by $S$.
When  $3g-3+s>0$, it carries a complete hyperbolic metric of finite volume, the so-called Poincar\'e metric; the points $p_1,
p_2, \cdots, p_s$ then become cusps at infinity. Even in the remaining cases, that is, for a once or twice punctured sphere, we can equip $S$ with a metric that is hyperbolic in the vicinity of the cusp(s), and for our purposes, the behavior of the metric there is all what counts, and we call such a metric Poincar\'e-like. In any case, our metric on $S$ is denoted by $\omega$.  Denote the
inclusion map of $S$ in $\overline S$ by $j$. Let $\rho: \pi_1(S)\to
Sl(n, \mathbb{C})$ be a  semisimple linear representation of
$\pi_1(S)$ which is unipotent near the cusps (for the precise
definition, cf. \S2.1). Corresponding to such a representation
$\rho$, one has a local system $L_{\rho}$ over $S$ and a
$\rho$-equivariant harmonic map $h: S\to Sl(n, \mathbb{C})/SU(n)$
with a certain special growth condition near the divisor. For the present  case of complex dimension 1, this is elementary; it also follows from the general result of \cite{jz-96}, see also  the remark in
\S2.2). This harmonic map can be considered as a Hermitian
metric on $L_{\rho}$---harmonic metric---so that we have a
so-called harmonic bundle $(L_{\rho}, h)$ \cite{sim-90}. Such a
bundle carries interesting structures, e.g. a Higgs bundle structure
$(E, \theta)$, where $\theta=\partial h$, and it has a
$\log$-singularity at the divisor.

The purpose of this note is to investigate various cohomologies of
$\overline S$ with degenerating coefficients $L_{\rho}$ (considered as a local system --- a flat vector bundle, a Higgs
bundle, or a $\mathcal{D}$-module, depending on the context): {\bf the \v{C}ech cohomology of $j_*L_{\rho}$}
(note that in the higher dimensional case, one needs to consider
the corresponding intersection cohomology \cite{cgm}), {\bf the
$L^2$-cohomology, the $L^2$-Dolbeault cohomology, and the
$L^2$-Higgs cohomology, and the relationships between them}. Here,
$L^2$ is defined by using the Poincar\'e(-like) metric $\omega$ and
the harmonic metric $h$. We want to generalize the results
\cite{zuc-79} valid for the case of variations of Hodge structures (VHS)
to the case of harmonic bundles, as was suggested by Simpson
\cite{sim-90}; in principle, in view of our assumption on the  representations in
question being unipotent, the situation should be similar to the
case of VHS.

This paper is meant to be a part of the general program of studying cohomologies with degenerating coefficients on quasiprojective varieties and their K\"ahlerian generalizations. The general aim here is
 not  restricted to the case
of curves nor to the one of representations that
are unipotent near the divisor. The purpose of this note therefore is to illuminate at this particular case where many of the (analytic and geometric) difficulties of the general case are not present what differences will appear when we consider
unipotent harmonic bundles instead of VHSs; for the case
of VHSs, the  various cohomologies have been considered by
various authors \cite{cks, lo, ss, yan1} and are well understood by now.

\vskip .3cm\noindent{\it Acknowledgements:} This work was begun
when the second author was visiting Mainz University, and its main part was
finished when he stayed at the Max Planck Institute for
Mathematics in the Sciences in Leipzig. He thanks  the above Institutes
for hospitality and a good working environment.

\section{The geometry associated with representations of fundamental groups}
\subsection{The decomposition of a flat connection}
In order to make this note an introduction into our general program, we describe some background material here. The knowledgable reader may skip this and perhaps also the next \S\S. Until the noncompact case is addressed explicitly below, the space $X$ will be compact in order to avoid difficulties with the analysis, in particular with the existence of the harmonic map.

Let $X$ be a Riemannian manifold, $V$ a $Gl(n,\mathbb{C})$ bundle on $X$ with a flat connection $D$ or, equivalently, a representation
$$ \rho:\pi_1(X) \to Gl(n,\mathbb{C}).$$
A metric $\langle \cdot,\cdot \rangle$ on $V$ is equivalent to a $\rho$-equivariant map
$$ h:\tilde{X} \to Gl(n, \mathbb{C})/U(n)=:\tilde{Y}$$
$(\langle v,v \rangle =h_{ij}v^i v^j)$.

The simplest case is, of course, $n=1$. Then $V$ is a line bundle, and
$$Gl(1, \mathbb{C})/U(1)=\mathbb{C}^*/S^1\simeq \mathbb{R}^+.$$
Using the isomorphism
$$ \log:\mathbb{R}^+ \to \mathbb{R},$$
a metric can be written as
$$ h=e^\lambda, \text{ for } \lambda:\tilde{X} \to \mathbb{R},$$
and so $h^{-1}dh =d\lambda$.\footnote{$h^{-1}dh$ is of course the derivative of the map $h:\tilde{X} \to \tilde{Y}$; usually, we should write $dh$, but the problem is that $d$ has two meanings, namely on one hand, the exterior derivative $d$, leading to the notation $h^{-1}dh$, and on the other hand, the differential $d$ of a map between Riemannian manifolds, suggesting to write $dh$. In the sequel, we use the notation $h^{-1}dh$ only for the case $n=1$ and write $dh$ else.}

Given a(ny) metric $h$ on $V$, the flat connection $D$ will not preserve it in general, and so, we split $D$ as
$$D=D_h +\theta,$$
where $D_h$ preserves $h$, i.e.,
$$d\langle v,v\rangle = \langle D_h v,v\rangle +\langle  v,D_hv\rangle.$$
Thus, we have $\theta =0$ iff $D$ is unitary iff $h$ is constant. The derivative $dh$ thus measures the deviation of $D$ from being unitary, i.e., $dh=\theta$. Thus, the energy of the map $h$ is given by $\int \|\theta\|^2$, and $h$ is harmonic iff
\begin{equation}
\label{harm}
D_h^* \theta =0.
\end{equation}
For $n=1$, reverting to our $h^{-1}dh$ notation, we have
$$D_h=d-h^{-1}dh, \text{ i.e. } \theta = h^{-1}dh=d\lambda.$$
The harmonic map equation (\ref{harm}) becomes
$$ 0=(d^*-h^{-1}dh)d\lambda=d^*d\lambda-d\lambda \wedge d\lambda = d^*d\lambda,$$
and so $\lambda$ is a harmonic function on $\tilde{X}$. $\lambda$ is not well defined as a function on $X$, but only as a map from $X$ into its Albanese variety. $d\lambda$, however, is well defined, and a harmonic 1-form, and $\lambda$ thus is the period map of a harmonic 1-form.

\subsection{Variations of Hodge structures}
Let $Z$ and $X$ be K\"ahler manifolds, with $\dim Z> \dim X$, $f:Z \to X$ holomorphic with smooth fibers. For the fiber $Z_x$ over $x\in X$, we have the Hodge decomposition
$$H^k(Z_x)=\bigoplus_{p+q=k}H^{p,q}(Z_x).$$
This induces a filtration ($n=\dim_{\mathbb{C}} Z_x=\dim_{\mathbb{C}}Z- \dim_{\mathbb{C}}X$)
$$0\subset F_x^n \subset \dots \subset F_x^{q+1} \subset F_x^q\subset \dots \subset H^k(Z_x).$$
This filtration defines an element of a subdomain of a Grassmannian of flags. From this, one obtains the Griffiths period map
\begin{equation}
\label{per1}
X \to D
\end{equation}
into the Griffiths period domain that is obtained by also imposing a polarization, i.e., a Hermitian form $\langle u,v \rangle$ s.t.
\begin{equation}
\label{perp}
H^{p,q} \perp H^{r,s} \text{ for } (p,q)\ne (r,s)
\end{equation}
and
\begin{equation}
\label{indef}
 (-1)^p\langle v,v \rangle >0 \text{ for } v\in H^{p,q}.
\end{equation}
This Griffiths period domain admits a holomorphic map
\begin{equation}
\label{per2}
D\to G/K
\end{equation}
onto a Hermitian symmetric space of noncompact type. $D$ is itself not K\"ahler because the natural metric coming from the polarization is indefinite, see (\ref{indef}). However, the image of the period map is tangential to a holomorphic distribution that is not integrable, but has nonpositive curvature in the tangential directions (those correspond to the directions in $G/K$).

There is a natural flat connection $D$ on the above decomposition, obtained by translating cohomology classes topologically or, equivalently, by considering the flat vector bundle $V$ with fiber
$$V_x=H^k(Z_x)$$
over $x \in X$. If we relate $D$ to the above decomposition, we obtain
\begin{equation}
\label{decomp}
D=\partial + \bar{\partial}+\theta +\bar{\theta}:H^{p,q} \to \Omega^{1,0}(H^{p,q})\oplus\Omega^{0,1}(H^{p,q})\oplus\Omega^{1,0}(H^{p-1,q+1})\oplus\Omega^{0,1}(H^{p+1,q-1}).
\end{equation}
$\partial$ and $\bar{\partial}$ simply come from the complex structure. $\theta$ and $\bar{\theta}$ have to shift the degree because in contrast to $\partial$ and $\bar{\partial}$, they do not operate by differentiation, but rather by multiplication with a 1-form -- as always when we split a connection as $D=d+A$.

Now an abstract complex variation of Hodge structures (VHS) is defined as a complex vector bundle $V$ over $X$ with a decomposition
$$V=\bigoplus_{p+q=k}V^{p,q}$$
as in (\ref{perp}, \ref{indef}) and a flat connection $D$ satisfying (\ref{decomp}).

This leads to a holomorphic vector bundle
$$E=\bigoplus_{p+q=k}E^{p,q},$$
using the above $\bar{\partial}$, with an endomorphism valued 1-form
$$\theta:E^{p,q}\to E^{p-1,q+1}\otimes\Omega^1(X)$$
with
$$\theta \wedge \theta =0$$
(this is the flatness condition on $D$).

By a result of Griffiths, $E$ is stable, and the flatness implies $c_i(E)=0$ for all $i$.

\subsection{Harmonic bundles and Higgs bundles}
It was then Simpson's fundamental idea (\cite{sim-92}) to revert this construction. A Higgs bundle $(E,\theta)$ over $X$ consists of a holomorphic vector bundle $E$ and $\theta:E\to E\otimes\Omega^1(X)$ with $\theta \wedge \theta =0$. The understanding of the relationship between harmonic bundles and Higgs bundles also owes much to Hitchin's important paper \cite{hit-87}.

By a theorem of Narasimhan-Seshadri for curves and by Donaldson, Uhlenbeck-Yau, Simpson for higher dimensional $X$, for a stable Higgs bundle, one can construct a Hermitian Yang-Mills connection $D_0$, and $D_0+\theta$ is flat if all $c_i(E)=0$. Such a $D_0$ then defines a harmonic metric on $E$, i.e., a harmonic map into a symmetric space $G/K$ as in (\ref{per1}, \ref{per2}). Conversely, from a semisimple representation
$$\rho: \pi_1(X) \to G,$$
$G$ a linear algebraic group, one obtains a $\rho$-equivariant harmonic map
$$h:\tilde{X} \to G/K,$$
and this defines a Higgs bundle $(E,\theta)$ with $\theta=dh$ ($\theta \wedge \theta =0$ follows from the pluriharmonicity of $h$ as originally discovered in a somewhat different context by Jost-Yau \cite{jy-83}; in fact, in our case of a Riemann surface, $\theta \wedge \theta =0$ is trivial because $\theta$ is a $(1,0)$-form, see below). More precisely,
decompose our flat connection $D=d'+d''$ into operators of type $(1, 0)$ and
$(0,1)$ respectively. Let $\delta'$ and $\delta''$ be the unique
operators of type $(1, 0)$ and $(0,1)$ such that the connections
$\delta'+d''$ and $\delta''+d'$ preserve the metric $h$. Let
$\partial=(d'+\delta')/2$ and
$\overline\partial=(d''+\delta'')/2$, and let
$\theta=(d'-\delta')/2$ and $\overline\theta=(d''-\delta'')/2$. It
is clear that $\mu=\partial+\overline\partial$ also preserves the
metric and $\overline\theta$ is the conjugate adjoint of $\theta$
w.r.t. $h$. Also, as mentioned,
$(\overline\partial+\theta)^2=0$.
Thus we obtain a structure of Higgs bundle on the bundle $L_{\rho}$ defined by the representation $\rho$ with
$\overline\partial$ being the holomorphic structure and $\theta$
the Higgs field; later on we denote it by $(E,
\overline\partial+\theta)$.

We also have the {\bf K\"ahler identities:} Set
$D''_h=\overline\partial+\theta$, $D'_h=\partial+\overline\theta$
and $D_h^c=D''_h-D'_h$. Note that $D=D'+D''$ and
$D_h''=(D+D_h^c)/2$. Let $\Lambda$ be the adjoint of the operation
of wedging with the K\"ahler form $\omega$. Then one has the first
order K\"ahler identities
\begin{eqnarray*}
&&(D'_h)^*=\sqrt{-1}[\Lambda, D''_h],
~~~(D''_h)^*=-\sqrt{-1}[\Lambda,
D'_h]\\
&&(D^c_h)^*=-\sqrt{-1}[\Lambda, D], ~~~(D)^*=\sqrt{-1}[\Lambda,
D^c_h],
\end{eqnarray*}
where $^*$ represents the adjoint of the respective operator. Set
$\Delta=DD^*+D^*D$ and $\Delta''=D''_h(D''_h)^*+(D_h'')^*D''_h$.
Using the above first order identities, one then has
\[
\Delta=2\Delta''.
\]
This shows that spaces of $\Delta$-harmonic forms valued in the
local system $L_{\rho}$ can be identified with that of
$\Delta''$-harmonic forms valued in the Higgs bundle $E$.

\subsection{The noncompact case}
When $X$ is no longer compact, the geometry works as before, but there arise difficulties with  the existence of the harmonic map $h$. On the other hand, the geometry of the bundle near a compactifying divisor leads to very interesting structures which are, in fact, our main interest. We first turn to the analytic  aspect. Even though in the present paper we shall be mainly concerned with the case of curves where the existence result is essentially elementary, in order to put the paper into a proper perspective, we shall describe the equivariant harmonic maps
(equivalently harmonic metrics) due to Jost-Zuo \cite{jz-96} when
the representations of $\pi_1$ are linearly semisimple and
unipotent and some consequences, all of which will be used in the
next section. The construction of the most general equivariant maps
corresponding to general linearly reductive representations will
be given in \cite{yan}; especially, there we will see that the
construction of Jost-Zuo \cite{jz-96, jz-97} corresponds to taking
the trivial filtration structure on the corresponding local
systems (cf. \cite{sim-90}).

Throughout this note, we will take the Poincar\'e-like metric on
the base manifolds, namely on punctured disks $\Delta^*$ near the
divisor (puncture) the metric is isometric to
\[
{\frac{dz\wedge d\overline{z}}{|z|^2(\log|z|)^2}}.
\]
Such a metric is complete, of finite volume and bounded geometry.

Let $\rho: \pi_1(S)\to {GL}(n, \mathbb{C})$ be a
semisimple linear representation, and restrict $\rho$ to a neighborhood
of $p_i\in D$, say a disk around $p_i$, which we call the {\it
boundary representation} of $\rho$. Throughout this note, {\it we
assume that all such boundary representations of $\rho$ are
unipotent}. This means that if denoting the image under $\rho$ of
the generator $\pi_1(\Delta^*)$ by the matrix $\gamma$, then
$\gamma$, under a suitable basis, can be represented by a
upper-triangle matrix with all the diagonal entries being $1$. It
is worth pointing out that in \cite{jz-96} the authors used a
 geometric definition  since their construction needs to
apply to more general target spaces; if the representations of
$\pi_1$ are linear, it is easy to see that their geometric
condition  reduces to our unipotent condition. In any case, from a geometric perspective, the essence of this condition is the following. Let $c$ be a closed curve in $\Delta^*$ representing $\rho$, for example the circle $r=r_0$, $r$ again being the Euclidean radius in $\Delta^*=\{ z\in \mathbb{C}:0<|z|<1\}$ (0 being the puncture, cusp, divisor,...). Thus $c$ is freely homotopic to arbitrarily short curves as we can move it closer to 0. In geometric terms, the unipotency condition means that the image $\gamma$ of our homotopy class can also be represented by arbitrarily short curves. The prototype is the representation $\rho$ corresponding to the identity map of $\Delta^*$, equipped with the Poincar\'e metric $\omega$. Since that metric has finite volume, the identity map has finite energy. Conversely, when the lengths of all curves homotopic to the image of $c$ have a positive lower bound, then any map in that homotopy class has infinite energy. An example are maps $f:\Delta^* \to U(1)=S^1$ equivariant w.r.t. the representation $\rho:\pi_1(\Delta^*)\to \mathbb{Z}$ that maps the generator of $\pi_1(\Delta^*)$ to the generator of $\pi_1(U(1))=\pi_1(S^1)=\mathbb{Z}$. In that case, any map homotopic to $f$ necessarily has infinite energy. In particular, while one can also show the existence of a harmonic map in this -- non-unipotent -- case, that map does not have finite energy, and the analysis becomes more subtle.
\begin{prop}
Let $\rho: \pi_1(S)\to Sl(n, \mathbb{C})$ be a semisimple
representation all boundary representations of which are
unipotent. Then there exists a $\rho$-equivariant harmonic map of
finite energy
\[
h: {\tilde S}\to Sl(n, \mathbb{C})/SU(n),
\]
where $\tilde S$ is the universal covering of $S$; moreover the
norm of the derivative $dh$ of $h$, when going down to $S$ and
measured near the divisor with respect to the Poincar\'e-like
metric and the standard Riemannian symmetric metric on $Sl(n,
\mathbb{C})/SU(n)$, $\le C{|\log r|^2}$ for some constant $C>0$,
where $r$ is the Euclidean radius of $\Delta^*$. Therefore, in
Simpson's notation \cite{sim-90}, the harmonic bundle $(L_{\rho},
h)$ is tame.
\end{prop}

Let $(L_{\rho}, h)$ be the harmonic bundle in the Proposition 1
and $D$ the flat connection. Since the harmonic bundle is tame, so
our discussions in the following and the next subsection lie in
the framework of Simpson \cite{sim-90}. We first describe the  norm
estimate  for flat sections near the divisor w.r.t. the
harmonic metric $h$ (for details, cf. \cite{sim-90}). We restrict
ourself to the punctured disk $\Delta^*$. Denote the image of the
generator of $\pi_1(\Delta^*)$ by $\gamma$, called the monodromy
of $L_{\rho}$, which is by the assumption unipotent; denote its
logarithm by $N$---the logarithmic monodromy, which is
nilpotent. Canonically, each fiber of $L_{\rho}$ has a so-called
weight filtration $\{W_l\}_{l=-k}^k$ arising from $N$ ($k$ is the
weight of $N$), which is $D$-invariant and which therefore
determines a filtration of $L_{\rho}$ by some local subsystems,
denoted by ${\bf W}_l$. Then a section $v$ of ${\bf W}_l$, if not
lying in ${\bf W}_{l-1}$, has the following norm estimate
\[
\Vert v\Vert^2_h\sim|\log r|^l.
\]

As explained above, canonically, from the
harmonic bundle $(L_{\rho}, h)$  constructed in the Proposition (using the estimate for
the derivative of $h$), one can derive a structure of Higgs bundle
on $L_{\rho}$ \cite{jz-97}.

Again since $(E, \overline\partial+\theta)$ comes from a tame
harmonic bundle, Theorem 1 of \cite{sim-90} works; in particular,
the curvature $R_{\mu}$ corresponding to the connection
$\mu=\partial+\overline\partial$ is bounded and the Higgs field
$\theta$ has a $\log$-singularity, and hence by a result of
Griffiths-Cornalba \cite{cg} one can extend $E$ across the
divisor, denoted by $j_*E$ as usual. More importantly, we have an
identification between the weight filtration of $N$ (the
logarithmic monodromy of $L_{\rho}$) on $L_{\rho}$ and the weight
filtration $\{W_l\}_{l=-k}^k$ of the residue of $\theta$ on
$res_{p_i}(j_*E)$ (Simpson's notation); this derives a similar norm
estimate  for meromorphic sections of $E$, i.e. let $e$ be
a meromorphic section, if the image in $res_{p_i}(j_*E)$ of $e$
lies in $W_l$ but not in $W_{l-1}$, then one has
\[
\Vert e\Vert^2\le C|\log r|^l.
\]
{\it Remark.} Due to the construction of the harmonic maps of
\cite{jz-96, jz-97} and the assumption that the representation is
unipotent, we have only the trivial filtered structure
\cite{sim-90} at the divisor on both the local system $L_{\rho}$ and
the Higgs bundle $(E, D'')$. And hence in the norm estimates
for flat sections of $L_{\rho}$ we do not have factors of the form
$r^b, b\in\mathbb{R}$; similarly for the estimates of meromorphic
sections of $(E, D'')$ one does not  have factors like
$r^{\alpha}, 0<\alpha< 1$. (For details, cf. \cite{sim-90, yan1})

\vskip .2cm\noindent
\vskip .2cm\noindent {\it Remark.} It should be pointed out that
in order that  various dualities and integration make sense, one
needs to restrict various discussions above to related
$L^2$-spaces, as is done in the next section, since our base
manifolds are noncompact; otherwise we will not have nice
identities. We omit these details, since they are standard
by now.

\section{The $L^2$-cohomology and the $L^2$-Higgs cohomology}
In this section, we continue to assume that the representation
$\rho$ is semisimple and unipotent. Therefore we have the
corresponding harmonic bundle $(L_{\rho}, h)$ and the Higgs bundle
$(E, \theta)$ with $\theta=\partial h$ and the hermitain metric
$h$; we also have various norm estimates   w.r.t. the metric $h$;
the most important thing is that one can identify the logarithmic
monodromy and the residue of $\theta$.

\subsection{The $L^2$-cohomology: The $L^2$-Poincar\'e Lemma}
Denoting the inclusion map of $S$ in $\overline S$ by $j$, one has
the direct image sheaf $j_*L_{\rho}$ of the local system
$L_{\rho}$ on $\overline S$ and then the \v{C}ech cohomology
$H^*({\overline S}, j_*L_{\rho})$. On the other hand, using the
Poincar\'e-like metric $\omega$ on $S$ and the harmonic metric $h$
on $L_{\rho}$, one can define a complex
$\{\mathcal{A}^{\cdot}_{(2)}(L_{\rho}), D\}$ of fine sheaves on
$\overline S$ as follows. Let $U$ be an open subset of $\overline
S$. Then $\mathcal{A}^{i}_{(2)}(L_{\rho})(U)$ is defined as the
set of

\vskip .3cm\noindent{\it $L_{\rho}$-valued $i$-forms $\eta$ on
$U\cap S$, with measurable coefficients and measurable exterior
derivative $D\eta$, such that $\eta$ and $D\eta$ have finite $L^2$
norm on $K\cap S$, for any compact subset $K$ of $U$, where $D$ is
the canonical flat connection of $L_{\rho}$.}

\vskip .3cm\noindent Since the sheaves are fine, so the cohomology
of the complex of global sections computes the hypercohomology of
$\{\mathcal{A}^{\cdot}_{(2)}(L_{\rho}), D\}$, namely
\[
\mathbb{H}^*({\overline S}, \{\mathcal{A}^{\cdot}_{(2)}(L_{\rho}),
D\})\cong H^*(\{\Gamma(\mathcal{A}^{\cdot}_{(2)}(L_{\rho})), D\}).
\]
We call this cohomology the {\it $L^2$-cohomology} of $\overline
S$ with values in $L_{\rho}$, denoted by $H^*_{(2)}({\overline S},
L_{\rho})$. The purpose of this subsection is then to establish
the following identification

\begin{thm}
There exists a natural identification
\[
H^*({\overline S}, j_*L_{\rho})\cong H^*_{(2)}({\overline S},
L_{\rho}).
\]
\end{thm}
{\it Remarks.}~ If $L_{\rho}$ comes from a variation of Hodge
structure (VHS), the identification was proved by S. Zucker
\cite{zuc-79}. In the higher dimensional case, instead of the
\v{C}ech cohomology, one needs to consider the intersection
cohomology \cite{cgm}; the identification was proved by
Cattani-Kaplan-Schmid \cite{cks} in the case of VHS. In general,
one has the following

\begin{conj}
There exists a natural identification
\[
H^*_{int}({\overline S}, j_*L_{\rho})\cong H^*_{(2)}({\overline
S}, L_{\rho}).
\]
\end{conj}

Canonically, the proof of Theorem 1 is reduced to prove

\begin{thm} (The $L^2$-Poincar\'e lemma)
The complex $\{\mathcal{A}^{\cdot}_{(2)}(L_{\rho}), D\}$ is a
resolution of $j_*L_{\rho}$. This is equivalent to saying that

\noindent 1)
$j_*L_{\rho}=\{\eta\in\mathcal{A}^{0}_{(2)}(L_{\rho})~|~D\eta=0\}$;

\noindent 2) the differential $D$ satisfies the Poincar\'e lemma,
i.e., if an $i$-form $\eta$ in $\mathcal{A}^{i}_{(2)}(L_{\rho})$
is $D$-closed, then there exists an $i-1$-form $\sigma\in
\mathcal{A}^{i-1}_{(2)}(L_{\rho})$ satisfying $D\sigma = \eta$,
for $i=1, 2$.
\end{thm}

In order to prove the above theorem, we need the following
\begin{lem}
Let $V$ be a constant one dimensional local system over
$\Delta^*$, with generator $v$, and assume that the corresponding
line bundle has a Hermitian metric with $\parallel
v\parallel^2\sim |\log r|^k$, where $r$ is the Euclidean radius.
Then the cohomology sheaves for $\mathcal{A}^{\cdot}_{(2)}(V)$
have stalks at the origin,
\begin{eqnarray*}
&&{\mathcal{H}}^0\big(\mathcal{A}^{\cdot}_{(2)}(V)\big) =
\begin{cases}
V ~~{\text{if}}~k\le 0 \\
0 ~~~{\text{if}}~k>0, \\
\end{cases}\\
&&{\mathcal{H}}^1\big(\mathcal{A}^{\cdot}_{(2)}(V)\big) =
\begin{cases}
{\frac{dt}{t}}\otimes V ~~~~~~~{\text{if}}~k\le -2 \\
0 ~~~~~~~~~~~~~~{\text{if}}~k\ge -1, k\neq 1\\
\mathcal{M}_1dr\otimes v ~~~{\text{if}}~k=1, \\
\end{cases}\\
&&{\mathcal{H}}^2\big(\mathcal{A}^{\cdot}_{(2)}(V)\big) =
\begin{cases}
0 ~~~~~~~~~~~~~~~~~~~~{\text{if}}~k\neq -1 \\
\mathcal{M}_1dr\wedge{\frac{dt}{t}}\otimes v ~~~{\text{if}}~k=-1, \\
\end{cases}
\end{eqnarray*}
where $\mathcal{M}_1$ is defined as
\[
{\frac{\{{\text{measurable functions}}~ f: \int^A_0|f(r)|^2|\log
r|(rdr)<\infty ~{\text{for some}}~ A<1\}}{\{f: f=u'~{\text{weakly
with}}~\int^A_0|u(r)|^2|\log r|^{-1}(r^{-1}dr)<\infty ~{\text{for
some}}~ A<1\}}}.
\]
\end{lem}
\noindent{\it Proof.} cf. \cite{zuc-79}, Proposition 6.6.

\vskip .3cm \noindent{\it Proof of Theorem 2.} On $S$, the
exactness is canonical. So the trouble comes from the singular
points $p_1, \cdots, p_s$, and since we are working with sheaves,
 from now on we just localize the problem to a punctured disk
$\Delta^*\subset\Delta$.

\vskip .3cm {\it The proof of 1).} This is equivalent to showing
that an $L^2$ flat section should be a section of
$j_*L_{\rho}$---an invariant section of $\gamma$, which
equivalently lies in the kernel of $N$.
To this end, denote the image of the generator of
$\pi_1(\Delta^*)$ under the representation $\rho$ by $\gamma$,
which, by the assumption, is unipotent, $\log\gamma$ by $N$, which
is nilpotent. For a (multi-valued) flat section $v$ of $L_{\rho}$,
setting
\[
{\tilde v}=\exp({\frac{1}{2\pi\sqrt{-1}}}N\log t)v,
\]
which is single-valued and $d''$-holomorphic ($D=d'+d''$, the $(1,
0)$ and $(0, 1)$-part respectively), one then has the canonical
extension $\overline{L}_{\rho}$ of $L_{\rho}$ to ${\overline S}$
($=\Delta$ locally) when $L_{\rho}$ is considered as a
$d''$-holomorphic bundle: sections of $\overline{L}_{\rho}$ at the
origin are generated by sections of the form $\tilde v$. Since $N$
is nilpotent, each fiber of $L_{\rho}$ canonically has a so-called
weight filtration $\{W_i\}$ which is $D$-invariant and which
therefore determines a filtration of $L_{\rho}$ by some local
subsystems, denoted by ${\bf W}_i$ and the corresponding extension
by $\overline{\bf W}_i$.

By the norm estimate of the harmonic metric $h$, it is easy to see
that a $d''$-holomorphic $L^2$ section of $L_{\rho}$ on $\Delta$
lies in $\overline{\bf W}_0+t\overline{L}_{\rho}$. On the other
hand, it is clear that $d'$-closed sections of $\overline{\bf
W}_0+t\overline{L}_{\rho}$ should be generated by sections of the
form $\tilde v$ satisfying $Nv=0$. This finishes the proof of 1).

\vskip .3cm {\it The proof of 2).} Without loss of generality, we
can assume that the representation $\rho: \pi_1(\Delta^*)\to Gl(m+1,
\mathbb{C})$ is irreducible; equivalently, $N$ acts irreducibly on
$\mathbb{C}^{m+1}$.
Canonically, one has a weight filtration of $\mathbb{C}^{m+1}$:
\[
0\subset W_{-m}\subset W_{-(m-2)}\subset\cdots\subset
W_{m-2}\subset W_{m}=\mathbb{C}^{m+1},
\]
which satisfies that the quotient ${\text{Gr}}^W_{i}=W_i/W_{i-2}$
is of dimension $1$ and $NW_i=W_{i-2}$. Correspondingly, one has
the invariant subbundles of $L_{\rho}$,
\[
0\subset {\bf W}_{-m}\subset {\bf W}_{-(m-2)}\subset\cdots\subset
{\bf W}_{m-2}\subset {\bf W}_{m}=L_{\rho},
\]
and hence the corrsponding filtration of
$\mathcal{A}^{\cdot}_{(2)}(L_{\rho})$
\[
0\subset\mathcal{A}^{\cdot}_{(2)}({\bf
W}_{-m})\subset\mathcal{A}^{\cdot}_{(2)}({\bf
W}_{-(m-2)})\subset\cdots\subset\mathcal{A}^{\cdot}_{(2)}({\bf
W}_{m-2})\subset\mathcal{A}^{\cdot}_{(2)}(L_{\rho}),
\]
which is clearly $D$-invariant. For simplicity, we denote
$\mathcal{A}^{\cdot}_{(2)}({\bf W}_{-i})$ by $K_i$; one has then a
filtered complex
$$
K_{-m}\supset K_{-m+2}\supset\cdots\supset K_{m-2}\supset
K_m\supset 0;
$$
note that $K_{-m+1}=K_{-m+2}, \dots$. Consider the quotient
$K_i/K_{i+1}$, which is clearly
$\mathcal{A}^{\cdot}_{(2)}({\text{Gr}}^W_{-i}(L_{\rho}))$. We now
have the spectral sequence $(E_r, d_r)_{r\ge 1}$ of the filtered
complex $\{K_i, D\}_{i=-m}^{m}$, which, by the theory of spectral
sequences, converges to the cohomology of $K_{-m}$, namely the
sheaf cohomology of $\{\mathcal{A}^{\cdot}_{(2)}(L_{\rho}), D\}$.
To prove Theorem 2, we need to analyze the sequence $(E_r,
d_r)_{r\ge 1}$.

By the definition of spectral sequences,
\[
E_1^{p,q}=H^{p+q}({\text{Gr}}^pK^*)=
\mathcal{H}^{p+q}\big(\mathcal{A}^{\cdot}_{(2)}({\text{Gr}}^W_{-p}(L_{\rho}))\big),
\]
where ${\text{Gr}}^pK_*=K_p/K_{p+1}$. We will show that the
differential $d_1: E^{p,q}_1\to E^{p+1,q}_1$ is trivial, as can be
observed as follows. Observe the diagram
\begin{equation*}
\begin{CD}
0\longrightarrow\mathcal{A}^{p+q}_{(2)}({\bf W}_{-p-1})@>i>>
\mathcal{A}^{p+q}_{(2)}({\bf W}_{-p})@>{Proj}>>
\mathcal{A}^{p+q}_{(2)}({\text{Gr}}^W_{-p}(L_{\rho}))@>>>{0}\\
@VV{D}V @VV{D}V @VV{D}V \\
0\longrightarrow\mathcal{A}^{p+q+1}_{(2)}({\bf W}_{-p-1})@>i>>
\mathcal{A}^{p+q+1}_{(2)}({\bf W}_{-p})@>{Proj}>>
\mathcal{A}^{p+q+1}_{(2)}({\text{Gr}}^W_{-p}(L_{\rho}))@>>>{0}\\
\end{CD}
\end{equation*}
where $i$ and $Proj$ are the inclusion and the projection
respectively; the third $D$ is the derived one. Let
$\phi\otimes\tilde v\in\mathcal{A}^{p+q}_{(2)}({\bf W}_{-p})$
represent a cohomologicy class from
$\mathcal{H}^{p+q}\big(\mathcal{A}^{\cdot}_{(2)}({\text{Gr}}^W_{-p}(L_{\rho}))\big)$;
equivalently, $d\phi=0$ (since $D\tilde v\in\Omega^1({\bf
W}_{-p-1})$, more precisely $\in\Omega^1({\bf W}_{-p-2})$). So
$D(\phi\otimes\tilde v)\in\mathcal{A}^{p+q+1}_{(2)}({\bf
W}_{-p-2})$, which represents a trivial cohomology class in
$\mathcal{H}^{p+q+1}\big(\mathcal{A}^{\cdot}_{(2)}({\text{Gr}}^W_{-p-1}(L_{\rho}))\big)(=E^{p+1,
q}_1)$, namely $d_1([\phi\otimes\tilde v])=0$.

Next, let us consider $d_2$. Since $d_1=0$, $E^{p,q}_2=E^{p,q}_1$.
Similar to the above argument, taking a form $\phi\otimes\tilde
v\in\mathcal{A}^{p+q}_{(2)}({\bf W}_{-p})$, which represents a
cohomology class of
$\mathcal{H}^{p+q}\big(\mathcal{A}^{\cdot}_{(2)}({\text{Gr}}^W_{-p}(L_{\rho}))\big)$,
$D(\phi\otimes\tilde
v)=(-1)^{p+q}{\frac{1}{2\pi\sqrt{-1}}}\big((\phi\wedge{\frac{dt}{t}})\otimes
N\tilde v\big)$ and then lies in $\mathcal{A}^{p+q+1}_{(2)}({\bf
W}_{-p-2})$, which represents a cohomological class of
$E^{p+2,q-1}_2=\mathcal{H}^{p+q+1}\big(\mathcal{A}^{\cdot}_{(2)}({\text{Gr}}^W_{-p-2}(L_{\rho}))\big)$.
So, $d_2$ is induced by
${\frac{1}{2\pi\sqrt{-1}}}{\frac{dt}{t}}N$.

Let us now consider the kernel and image of $d_2$. Applying Lemma 1
to
$E^{p,q}_2=\mathcal{H}^{p+q}\big(\mathcal{A}^{\cdot}_{(2)}({\text{Gr}}^W_{-p}(L_{\rho}))\big)$,
we get that  the only (possibly) nontrivial terms at $E_2$ are
$\{E^{p,-p}_2\}_{p\ge 0}$, $\{E^{p+2, -p-1}_2\}_{p\ge 0}$,
$E^{-1,2}_2$, and $E^{1, 1}_2$. Furthermore, from the above
argument together with Lemma 1, we obtain that $d_2: E^{p,-p}_2\to
E^{p+2, -p-1}_2, {p\ge 0}$ (if $E^{p+2, -p-1}_2$ is nontrivial,
i.e. $p\le m-2$) and $d_2: E^{-1,2}_2\to E^{1, 1}_2$ are
isomorphisms and that $d(E^{p+2, -p-1}_2)=0, p\ge 0$ and
$d_2(E^{1,1}_2)=0$.

Summing all the above up,  the only possible
nontrivial terms at $E_3$ are $E^{m-1, -m+1}_3, E^{m,-m}_3$. Thus
the spectral sequence $\{E_r, d_r\}_{r\ge 1}$ degenerates at $E_3$
and  the only possible terms are $E^{m-1, -m+1}_3, E^{m,-m}_3$.
Therefore, by the theory of spectral sequences of filtered
complexes, $\mathcal{H}^i(\mathcal{A}^{\cdot}_{(2)}(L_{\rho}))=0,
i=1, 2$. The proof of 2) is finished.
~~~~~~~~~~~~~~~~~~~~~~~~~~~~~~~~~~~~~~~~~~~~~~~~~~~~~~~~~~~~~~~~~$\Box$

\subsection{The $L^2$-Higgs cohomology: The $L^2$-$\overline\partial$-Poincar\'e Lemma}
As seen in \S 2, the harmonic metric $h$ (the $\rho$-equivariant
harmonic map) with the tame growth condition on $L_{\rho}$ derives
a structure of Higgs bundle on $L_{\rho}$: $(E,
D''=\overline\partial+\theta)$, satisfying $D=D'+D''$ with
$D'=\partial+\overline\theta$; moreover $E$, as a holomorphic
vector bundle, has bounded curvature under the harmonic metric $h$
so that $E$ can be analytically extended to $\overline S$, as
usual just denoted by $j_*E$. Furthermore, $\theta$ has a
$\log$-singularity, i.e. $\theta\sim {\frac{dt}{t}}N$. It is
especially worth pointing out that by an argument of Simpson (cf.
\cite{sim-90}), {\it the residue $N$ of $\theta$ here coincides
with the logarithmic monodromy $N$ in the local system
$L_{\rho}$}; so although under different bundle structures, we
have the same weight filtration under certain suitable
identification. Throughout this subsection, we consider the Higgs
bundle $(E, D''=\overline\partial+\theta)$ together with the
harmonic metric $h$, just forgetting that it comes from the local
system corresponding to the representation $\rho$. As in the
previous subsection, using the Poincar\'e-like metric $\omega$ on
$S$ and the harmonic matric $h$ on $(E, D'')$, one can similarly
define a complex $\{\mathcal{A}^{\cdot}_{(2)}(E), D''\}$ of fine
sheaves on $\overline S$. Let $U$ be an open subset of $\overline
S$. Then $\mathcal{A}^{i}_{(2)}(E)(U)$ is defined as the set of

\vskip .3cm\noindent {\it $E$-valued $i$-forms $\eta$ on $U\cap
S$, with measurable coeffcients and measurable exterior derivative
$\overline\partial\eta$, such that $\eta$ and $D''\eta$ have
finite $L^2$ norm on $K\cap S$, for any compact subset $K$ of
$U$.}

\vskip .3cm\noindent {\it Remarks.} 1) The Higgs condition
$(D'')^2=0$ makes $\{\mathcal{A}^{\cdot}_{(2)}(E), D''\}$ a complex, which is actually a complex of certain
differential forms valued in $j_*E$. 2) Due to following lemma, it
is actually sufficient to assume that $\eta$ and
$\overline\partial\eta$ are $L^2$ in the above definition. 3)
Since the sheaves are fine, so again the hypercohomology
$\mathbb{H}^*(\{\mathcal{A}^{\cdot}_{(2)}(E), D''\})$ is computed
by the cohomology ${H}^*(\{\Gamma(\mathcal{A}^{\cdot}_{(2)}(E)),
D''\})$ of the complex of global sections; we call it the {\it
$L^2$-Higgs cohomology} of $\overline S$ valued in the Higgs
bundle $(E, D=\overline\partial+\theta)$, denoted by
$H^*_{(2)}({\overline S}, E)$.

\begin{lem}
$\theta$ is an $L^2$-bounded operator.
\end{lem}
{\it Proof}. As mentioned before, $\theta\sim{\frac{dt}{t}}N$. So
it suffices to show that ${\frac{dt}{t}}N$ is $L^2$-bounded. Since
$N$ lowers weights by $2$, so by the estimate of the norm of the
harmonic bundle, the norm changes by multiplication with $(\log|t|)^{-2}$;
on the other hand,
$\Vert{\frac{dt}{t}}\Vert^2_\omega=(\log|t|)^2$. The proof is
obtained. ~~~~~~~~~~~~~~~~~~~~~~~~~~~~~~~~~~~~~$\Box$

\vskip .3cm On the other hand, based on the above lemma one can
also define a sub-complex of $\{\mathcal{A}^{\cdot}_{(2)}(E),
D''\}$---the $L^2$-{\it holomorphic Dolbeault complex}
$\{\Omega^{\cdot}_{(2)}(E), \theta\}$ as follows:
$\Omega^{i}_{(2)}(E)(U)$ is defined as the set of

\vskip .3cm\noindent{\it $E$-valued holomorphic $i$-forms $\eta$
on $U\cap S$ such that $\eta$ has finite $L^2$ norm on $K\cap S$,
for any compact subset $K$ of $U$.}

\vskip .3cm\noindent $\theta\wedge\theta=0$ makes
$\{\Omega^{\cdot}_{(2)}(E), \theta\}$  again  a complex,
which is actually a complex of certain meromorphic differential
forms valued in $j_*E$. We call the hypercohomology
$\mathbb{H}^*(\{\Omega^{\cdot}_{(2)}(E), \theta\})$  the
$L^2$-{\it Dolbeault cohomology} of $\overline S$ valued in the
Higgs bundle $(E, D=\overline\partial+\theta)$.


The purpose of this subsection is then to show that the above two
complexes have the same hypercohomologies; more precisely
\begin{thm} ($L^2$-$\overline\partial$-Poincar\'e lemma)
The inclusion $i: \{\Omega^{\cdot}_{(2)}(E), \theta\}\hookrightarrow
\{\mathcal{A}^{\cdot}_{(2)}(E), D''\}$ is a quasi-isomorphism;
and hence one has
\[
\mathbb{H}^*(\{\Omega^{\cdot}_{(2)}(E),
\theta\})\cong{H}_{(2)}^*({\overline S}, E).
\]
\end{thm}
{\it Remarks.} In the case when $E$ comes from a VHS, the theorem was
showed by S. Zucker \cite{zuc-79} (for the case of curves) and
Jost-Yang-Zuo \cite{yan1} (for the general case).

\vskip .3cm In order to prove the theorem, we need some
preliminaries. First we give the following
\begin{defi}
Let $V$ be a Hermitian vector bundle over a Riemannian manifold
$M$, ${\bf v}=\{v_1, v_2, \cdots, v_q\}$ be a global frame field
of $V$. Then $\bf v$ is said to be $L^2$-adapted if
$\sum_{i=1}^qf_iv_i$ is square integrable implies that each
$f_iv_i$ is square integrable, where the $f_i$ are smooth functions
on $M$.
\end{defi}
We next show the following general lemma, which will be needed
when we consider general semisimple representations (not
necessarily being unipotent at infinity), although what we really
need here is the special case  $\alpha=0$ which was proved in
\cite{zuc-79}.
\begin{lem} Let $V$ be
a holomorphic line bundle on $\Delta^*$ with generating section
$\sigma$, and with a Hermitian metric satisfying
\[
\Vert\sigma\Vert^2\sim r^{\alpha}|\log r|^k, ~~0\le\alpha <1,
k\in\mathbb{Z}.
\]
Then for every germ of an $L^2$ $(0, 1)$-form
$\phi=fd\overline{t}\otimes\sigma$ at the origin, if 1) $\alpha=0$
and $k\neq 1$ or 2) $\alpha\neq 0$, there exists an $L^2$ section
$u\otimes\sigma$ with $\overline\partial u=fd\overline{t}$.
\end{lem}
{\it Proof.} Similar to \cite{zuc-79}, we also use Fourier series.
Using polar coordinates, we write $u$ and $f$ as $r$-dependent
Fourier series
\[
u=\sum_{n=-\infty}^\infty u_n(r)e^{\sqrt{-1}n\theta},
f=\sum_{n=-\infty}^\infty f_n(r)e^{\sqrt{-1}n\theta}.
\]
As $\partial/\partial\overline
t=(1/2)e^{\sqrt{-1}\theta}[\partial/\partial
r+(\sqrt{-1}/r)\partial/\partial\theta]$, the equation
$\overline\partial u=f$ becomes
\[
{\frac 1 2}\big(u'_n-{\frac n r}u_n\big)=f_{n+1}~~{\text{for
all}}~ n\in\mathbb{Z},
\]
for $C^\infty$ germs $u$ and $f$, or
\[
{\frac 1 2}{\frac{d}{dr}}[r^{-n}u_n(r)]=r_{-n}f_{n+1}(r).
\]
We are given that for some $A<1$,
\[
\Vert\phi\Vert^2_{(2)}=4\pi\sum_{n=-\infty}^{\infty}\int^A_0|f_n|^2|\log
r|^kr^{1+\alpha}dr<\infty,
\]
and we want to obtain some $u$ satisfying
\[
\Vert
u\otimes\sigma\Vert^2_{(2)}=2\pi\sum_{n=-\infty}^{\infty}\int^A_0|u_n|^2|\log
r|^{k-2}r^{-1+\alpha}dr<\infty.
\]
In order to obtain $u$, we make use of integral representations
and try to show that $\Vert u\Vert^2_{(2)}\le C\Vert
f\Vert^2_{(2)}$ for some positive constant $C$. Then a standard
approximation argument shows that the lemma is true. As pointed
out above, the case  $\alpha=0$ has been proved in
\cite{zuc-79}, so from now on we assume $\alpha>0$.  Set
\begin{eqnarray*}
u_n(r)=
\begin{cases}
1)~ 2r^n\int_0^r\rho^{-n}f_{n+1}(\rho)d\rho~~~~~~~~~{\text{if}}~n<0, \\
2)~ -2\int_r^Af_{1}(\rho)d\rho~~~~~~~~~~~~~~~{\text{if}}~n=0, \\
3)~ -2r^n\int_r^A\rho^{-n}f_{n+1}(\rho)d\rho~~~~~{\text{if}}~n>0. \\
\end{cases}
\end{eqnarray*}
(Note that the case of $\alpha>0$ is different from the case
$\alpha=0$. In the latter case, if $n=0$, one needs furthermore
consider the values of $k$; in the former case, we need not even
consider $k$, but one needs take $u_0=-2\int_r^Af_{1}(\rho)d\rho$,
since $\int_0^rf_{1}(\rho)d\rho$ can possible not be integrable.)
In order to prove $\Vert u\Vert^2_{(2)}\le C\Vert f\Vert^2_{(2)}$,
it is sufficient to prove that for all $n$
\[
\int_0^A|u_n|^2|\log r|^{k-2}(r^{-1+\alpha}dr)\le
C\int_0^A|f_{n+1}|^2|\log r|^{k}(r^{1+\alpha}dr).
\]
If $n<0$, one has
\begin{eqnarray*}
&&\int^A_0|u_n|^2r^{-1+\alpha}|\log r|^{k-2}dr\\
&\le& 4\int_0^Ar^{2n}\big(\int^r_0\rho^{-2n}|f_{n+1}(\rho)|^2d\rho
\big)\big(\int_0^rd\rho\big)r^{-1+\alpha}|\log r|^{k-2}dr\\
&=& 4\int_0^A\big(\int^r_0\rho^{-2n}|f_{n+1}(\rho)|^2d\rho
\big)r^{2n+\alpha}|\log r|^{k-2}dr\\
&=& 4\int_0^A\big(\int^A_\rho r^{2n+\alpha+1}|\log
r|^{k-2}(r^{-1}dr)\big)
\rho^{-2n}|f_{n+1}(\rho)|^2d\rho\\
&\le& 4\int_0^A\big(\int^A_\rho r^{-1}dr\big)
|f_{n+1}(\rho)|^2\rho^{\alpha+1}|\log
\rho|^{k-2}d\rho\\
&\le& C\int_0^A |f_{n+1}(\rho)|^2\rho^{\alpha+1}|\log
\rho|^{k}d\rho.
\end{eqnarray*}
In the second inequality above, we used the property of $|\log
r|^{k-2}r^{2n+\alpha+1}$ being decreasing since $2n+\alpha+1<0$.
If $n=0$, one has
\begin{eqnarray*}
|u_0(r)|^2&=&4\big(\int_r^A|f_1(\rho)|^2|\log\rho|^{1+k}\rho^{1+\alpha}
d\rho \big)\big(\int^A_r|\log\rho|^{-1-k}\rho^{-1-\alpha}dr\big)\\
&\le&C\big(\int_r^A|f_1(\rho)|^2|\log\rho|^{1+k}\rho^{1+\alpha}
d\rho \big)r^{-\alpha}|\log r|^{-k},
\end{eqnarray*}
in the above inequality, we again used the property of
$|\log\rho|^{-1-k}\rho^{-\alpha}$ being decreasing. Estimating
$\int_0^A|u_0(r)|^2|\log r|^{k-2}r^{-1+\alpha}dr$, we have
\begin{eqnarray*}
&&\int_0^A|u_0(r)|^2|\log r|^{k-2}r^{-1+\alpha}dr\\
&\le&C\int_0^A\big(\int_r^A|f_1(\rho)|^2|\log\rho|^{1+k}\rho^{1+\alpha}
d\rho \big)r^{-1}|\log r|^{-2}dr \\
&=&C\int_0^A\big(\int_0^r\rho^{-1}|\log\rho|^{-2}d\rho\big)|f_1(r)|^2|\log
r|^{1+k}r^{1+\alpha} dr\\
&=&C\int_0^A|f_1(r)|^2|\log r|^{k}r^{1+\alpha} dr.
\end{eqnarray*}
If $n>0$, one has
\begin{eqnarray*}
&&\int^A_0|u_n|^2r^{-1+\alpha}|\log r|^{k-2}dr\\
&\le&
4\int_0^Ar^{2n}\big(\int^A_r\rho^{-2n+1}|f_{n+1}(\rho)|^2d\rho
\big)\big(\int_r^A{\frac 1 \rho}d\rho\big)r^{-1+\alpha}|\log r|^{k-2}dr\\
&\le& C\int_0^A\big(\int^A_r\rho^{-2n+1}|f_{n+1}(\rho)|^2d\rho
\big)r^{2n+\alpha-1}|\log r|^{k-1}dr\\
&=& C\int_0^A\big(\int^\rho_0r^{2n+\alpha-1}|\log r|^{k-1}dr\big)
\rho^{-2n+1}|f_{n+1}(\rho)|^2d\rho\\
&\le&
C\int_0^A|f_{n+1}(\rho)|^2|\log\rho|^{k-1}\rho^{1+\alpha}d\rho.
\end{eqnarray*}
In the last inequality above, we use that
$r^{2n+\alpha-1}|\log r|^{k-1}$ is increasing, since
$2n+\alpha-1>0$. The proof of the lemma is completed. ~~~~~$\Box$

\vskip .3cm \noindent {\it Proof of Theorem 3.} Similar to the
proof of Theorem 2, the difficulty again comes from $\overline
S\setminus S$, so one can restrict the problem to a  small enough
punctured disk $\Delta^*$ so that we can choose a holomorphic
basis of $j_*E$: $e_1, e_2, \cdots, e_n$, which is compatible with
the weight filtation corresponding to $N$ (the residue of
$\theta$) at the origin. The compatibility implies that each $e_i$
of the basis has the property
\begin{eqnarray*}
\Vert e_i\Vert_h&\sim&|\log r|^k, ~{\text{for a certain integer}}~k;\\
\Vert \theta e_i\Vert_{h,\omega}&\sim&|\log r|^{k},
~{\text{if}}~\theta e_i\neq 0,
\end{eqnarray*}
where $\omega$ is the Poincar\'e-like metric. In particular, for
some section $e$, if $\Vert e\Vert_h\sim|\log r|$, then $\theta
e\neq 0$; if $\Vert e\Vert_h\sim|\log r|^{-1}$, then there exists
a section $e'$ satisfying $\theta e'={\frac{dt}{t}}\otimes e$.
Using these facts, it is not difficult to prove that the basis
$e_1, e_2, \dots, e_n$ is an $L^2$-adapted one.

In order to prove the theorem, we need to show that the inclusion
$i$ derives an isomorphism between the corresponding cohomology
sheaves at the origin; by the argument of standard homological
algebra, this is equivalent to showing that for any $D''$-closed form
$\phi\in{\mathcal{A}}^r_{(2)}(E)$ on a neighborhood $U$ of the
origin, there is a $\theta$-closed form $\eta\in\Omega^r_{(2)}(E)$
and a form $\psi\in{\mathcal{A}}^{r-1}_{(2)}(E)$ on (a possibly
smaller) $U$ satisfying $\phi=\eta+D''\psi$, $r=0, 1, 2$.

When $r=0$, it is clear that a $D''$-closed form is a holomorphic
section of $E$ and $\theta$-closed; in the case of $r=2$, the form
$\phi$ can be written as the sum of forms
$\phi'\wedge{\frac{dt}{t}}\otimes e$, $\phi'$ is a complex-value
$(0, 1)$-form and $e$ is an element of the basis. Since the basis
is $L^2$-adapted, each summand is still $L^2$. Considering
${\frac{dt}{t}}\otimes e$ as a generator of a holomorphic line
bundle, its norm  satisfies under the harmonic metric and
the Poincar\'e-like metric
\[
\Vert{\frac{dt}{t}}\otimes e\Vert=|\log r|^k,
\]
for some integer $k$. When $k\neq 1$, by the above lemma, there
exists an $L^2$-form $u{\frac{dt}{t}}\otimes e$ satisfying
\[
\overline\partial(u{\frac{dt}{t}}\otimes
e)=\phi'\wedge{\frac{dt}{t}}\otimes e,
\]
and hence $D''(u{\frac{dt}{t}}\otimes
e)=\phi'\otimes{\frac{dt}{t}}\otimes e$; assuming $k=1$, there
exists a holomorphic section $e'$ of $E$ satisfying $\theta
e'={\frac{dt}{t}}\otimes e$, and hence $\theta(-\phi'\otimes
e')=D''(-\phi'\otimes e')=\phi'\wedge{\frac{dt}{t}}\otimes e$. It
is clear that $-\phi'\otimes e'$ is $L^2$ (
$\overline\partial(-\phi'\otimes e')$ is automatically $L^2$ since
$=0$) and hence $-\phi'\otimes e'\in{\mathcal{A}}^{1}_{(2)}(E)$.
Summing all the above, we have that any form of
${\mathcal{A}}^{2}_{(2)}(E)$ is $D''$-coclosed.

We now turn to the case of $r=1$. In order to make the proof
clearer, we can refer from time to time to the following diagram
\begin{equation*}
\begin{CD}
\Omega^0_{(2)}(E)@>i>>
\mathcal{A}^{0,0}_{(2)}(E)@>{\overline\partial}>>
\mathcal{A}^{0,1}_{(2)}(E)@>>>{0}\\
@VV{\theta}V @VV{\theta}V @VV{\theta}V \\
\Omega^{1}_{(2)}(E)@>i>>
\mathcal{A}^{1,0}_{(2)}(E)@>{\overline\partial}>>
\mathcal{A}^{1,1}_{(2)}(E)@>>>{0}.\\
\end{CD}
\end{equation*}
Write $\phi=\phi^{1,0}+\phi^{0,1}$, $\phi^{1,0}$ (resp.
$\phi^{0,1}$) being the part of type $(1, 0)$ (resp. $(0, 1)$).
The $D''$-closedness of $\phi$ is equivalent to
$$
\overline\partial\phi^{1,0}+\theta\phi^{0,1}=0.
$$
Write $\phi^{0,1}$ as the sum of forms $\phi^{0,1}_i\otimes e_i$,
$i=1, \cdots, n$, which, by the $L^2$-adaptedness of the basis,
are $L^2$ again. Assume $\Vert e_{i_0}\Vert\sim|\log r|$ and for
$i\neq i_0$, $\Vert e_{i}\Vert\sim|\log r|^k, k\neq 1$. By the
above lemma, for $i\neq i_0$, there exists an $L^2$ section
$u_ie_i$ satisfying
$$
\overline\partial u_i\otimes e_i=\phi_i^{0,1}\otimes e_i.
$$
It is easy to see that the part of type $(0,1)$ of the form
$\phi-D''(u_ie_i)$ does not contain any further terms of the form
$\phi^{0,1}_i\otimes e_i, i\neq i_0$. So, w.l.o.g., we may assume
that $\phi^{0,1}=\phi_{i_0}^{0,1}\otimes e_{i_0}$.
In order to deal with the term $\phi_{i_0}^{0,1}\otimes e_{i_0}$,
we use the $D''$-closedness of $\phi$, i.e.
$\overline\partial\phi^{1,0}+\theta\phi^{0,1}=0$. Considering the
holomorphic vector bundle ${\frac{dt}{t}}\otimes E$, it is not
difficult to see that one can extend $\theta e_i, i=1, \cdots, n$,
to an $L^2$-adapted basis (we do not really need this basis). Due
to the relation $\overline\partial\phi^{1,0}+\theta\phi^{0,1}=0$,
one can write $\phi^{1,0}$ as $u_{i_0}\theta e_{i_0}$ satisfying
$\overline\partial u_{i_0}=\phi^{0,1}_{i_0}$. Clearly both
$u_{i_0}e_{i_0}$ and $\overline\partial(u_{i_0}e_{i_0})$ are
$L^2$, namely $u_{i_0}e_{i_0}\in{\mathcal{A}}^{0}_{(2)}(E)$, and
$\phi-D''(u_{i_0}e_{i_0})$ does no longer contain  the part of
type $(0, 1)$. Obviously a $D''$-closed $L^2$-form of type $(1,
0)$ is holomophic and $\theta$-closed. Summing the above argument
up, we have that for any $D''$-closed form
$\phi\in{\mathcal{A}}^1_{(2)}(E)$ on a neighborhood $U$ of the
origin, there is a $\theta$-closed form $\eta\in\Omega^1_{(2)}(E)$
and a section $\psi\in{\mathcal{A}}^{0}_{(2)}(E)$ on (a possibly
smaller) $U$ satisfying $\phi=\eta+D''\psi$. This finishes the
proof of the theorem.
~~~~~~~~~~~~~~~~~~~~~~~~~~~~~~~~~~~~~~~~~~~~~~~~~~~~~~~~~~~~~~~~~~~~~~~~~~~~~~~~~~~~~~~~~~~$\Box$

\vskip .3cm Using the K\"ahler identity for harmonic bundles (cf.
\S 2), $H^*_{(2)}(\overline S, L_{\rho})$ can be identified with
$H^*_{(2)}(\overline S, E)$, and hence we have the following
\begin{cor}
$$
H^*(\overline S,
j_*L_{\rho})\cong\mathbb{H}^*(\{\Omega^{\cdot}_{(2)}(E),
\theta\}).
$$
\end{cor}

\bigskip

\noindent Max-Planck Institute for Mathematics in the Sciences,
Leipzig \vskip .3cm\noindent
Department of Applied Mathematics, Tongji University, Shanghai\\
{\it E-mail}: {yhyang@tongji.edu.cn} \vskip .3cm\noindent
Department of Mathematics, Mainz University, Mainz
\end{document}